\documentclass[11pt]{amsart}
\usepackage{amsmath,amssymb,mathrsfs}
\newtheorem{theorem}{Theorem}
\newtheorem{proposition}[theorem]{Proposition}
\newtheorem{corollary}[theorem]{Corollary}
\newtheorem{lemma}[theorem]{Lemma}
\theoremstyle{definition}
\newtheorem{remark}[theorem]{Remark}

\begin{document}

\title[Quaternionic-K\"ahler manifolds]{Quaternionic-K\"ahler manifolds with nonnegative sectional curvature}
\author{Simon Brendle and Uwe Semmelmann}
\address{Columbia University \\ 2990 Broadway \\ New York NY 10027 \\ USA}
\address{Universit\"at Stuttgart \\ Pfaffenwaldring 57 \\ 70569 Stuttgart \\ Germany}
\thanks{The first author was supported by the National Science Foundation under grant DMS-2403981 and by the Simons Foundation. He acknowledges the hospitality of T\"ubingen University, where part of this work was carried out.}
\maketitle
\begin{abstract}
We show that a compact quaternionic-K\"ahler manifold with positive scalar curvature and nonnegative sectional curvature is isometric to a symmetric space. This extends a classical theorem of Berger.
\end{abstract}

\section{Introduction}

Throughout this paper, we assume that $(M,g)$ is a compact quaternionic-K\"ahler manifold of dimension $4m \geq 8$ with positive scalar curvature. In particular, there exists a parallel subbundle $\mathscr{G} \subset \text{\rm End}(TM)$ of rank $3$ with the following property. For each point $p \in M$, we can find linear transformations $I,J,K \in \mathscr{G}_p$ such that $I^2 = J^2 = K^2 = IJK = -\text{\rm id}$ and 
\[g(X,Y) = g(IX,IY) = g(JX,JY) = g(KX,KY)\] 
for all vectors $X,Y \in T_p M$. We may write 
\[\mathscr{G}_p = \{aI+bJ+cK: a,b,c \in \mathbb{R}\}.\] 
We define 
\[\mathscr{J}_p = \{aI+bJ+cK: \text{\rm $a,b,c \in \mathbb{R}$ and $a^2+b^2+c^2 = 1$}\} \subset \mathscr{G}_p.\] 
By work of Alekseevskii \cite{Alekseevskii} (see also \cite{Salamon}, Theorem 3.1), the Riemann curvature tensor of $(M,g)$ can be written in the form $R = R_1+\kappa R_0$. Here, $R_0$ is defined by 
\begin{align*} 
&4 \, R_0(X,Y,Z,W) \\ 
&= g(X,Z) \, g(Y,W) - g(X,W) \, g(Y,Z) \\ 
&+ 2 \, g(IX,Y) \, g(IZ,W) + g(IX,Z) \, g(IY,W) - g(IX,W) \, g(IY,Z) \\ 
&+ 2 \, g(JX,Y) \, g(JZ,W) + g(JX,Z) \, g(JY,W) - g(JX,W) \, g(JY,Z) \\ 
&+ 2 \, g(KX,Y) \, g(KZ,W) + g(KX,Z) \, g(KY,W) - g(KX,W) \, g(KY,Z) 
\end{align*}
for all vectors $X,Y,Z,W \in T_p M$ (cf. \cite{Besse}, equation 14.44). Moreover, $R_1$ is an algebraic curvature tensor satisfying 
\begin{align*} 
R_1(X,Y,Z,W) 
&= R_1(X,Y,IZ,IW) \\ 
&= R_1(X,Y,JZ,JW) \\ 
&= R_1(X,Y,KZ,KW) 
\end{align*} 
for all vectors $X,Y,Z,W \in T_p M$. To summarize, $R_1$ satisfies all the algebraic properties of a hyper-K\"ahler curvature tensor, including the first Bianchi identity. In particular, the Ricci tensor of $R_1$ vanishes identically. Consequently, the Ricci tensor of $(M,g)$ is given by $(m+2)\kappa \, g$.

By a theorem of Salamon, a compact quaternionic-K\"ahler manifold with positive scalar curvature is simply connected (see \cite{Salamon}, Theorem 6.6). Let $\text{\rm Hol}(M,g)$ denote the holonomy group of $(M,g)$. Since $M$ is simply connected, $\text{\rm Hol}(M,g)$ is a compact, connected Lie subgroup of $\text{\rm SO}(4m)$ (see \cite{Besse}, Corollary 10.47 and Theorem 10.51). Moreover, $\text{\rm Hol}(M,g) \subset \text{\rm Sp}(m) \cdot \text{\rm Sp}(1)$ since $(M,g)$ is a quaternionic-K\"ahler manifold. 

\begin{theorem}
\label{main.thm}
Let $(M,g)$ be a compact quaternionic-K\"ahler manifold of dimension $4m \geq 8$ with positive scalar curvature. We further assume that $R(X,Y,X,Y) + R(X,JY,X,JY) \geq 0$ for all points $p \in M$, all $I,J,K \in \mathscr{J}_p$ satisfying $I^2=J^2=K^2=IJK=-\text{\rm id}$, and all unit vectors $X,Y \in T_p M$ satisfying $g(X,Y)=g(IX,Y)=g(JX,Y)=g(KX,Y)=0$. If $\text{\rm Hol}(M,g) = \text{\rm Sp}(m) \cdot \text{\rm Sp}(1)$, then $R_1=0$ at each point in $M$. 
\end{theorem}

Combining Theorem \ref{main.thm} with Berger's holonomy classification theorem, we can draw the following conclusion.

\begin{corollary}
\label{consequence.of.main.thm}
Let $(M,g)$ be a compact quaternionic-K\"ahler manifold of dimension $4m \geq 8$ with positive scalar curvature. We further assume that $R(X,Y,X,Y) + R(X,JY,X,JY) \geq 0$ for all points $p \in M$, all $I,J,K \in \mathscr{J}_p$ satisfying $I^2=J^2=K^2=IJK=-\text{\rm id}$, and all unit vectors $X,Y \in T_p M$ satisfying $g(X,Y)=g(IX,Y)=g(JX,Y)=g(KX,Y)=0$. Then $(M,g)$ is isometric to a symmetric space.
\end{corollary}

Corollary \ref{consequence.of.main.thm} extends earlier work of Berger \cite{Berger2} and Chow and Yang \cite{Chow-Yang}. 

LeBrun and Salamon \cite{LeBrun-Salamon} have conjectured that a compact quaternionic-K\"ahler manifold with positive scalar curvature is isometric to a symmetric space. The LeBrun-Salamon conjecture is known to be true in dimension $8$ by work of Poon and Salamon \cite{Poon-Salamon}, and in dimension $12$ and $16$ by work of Buczy\'nski and Wi\'sniewski \cite{Buczynski-Wisniewski}.

\section{Proof of Theorem \ref{main.thm}}

In this section, we give the proof of Theorem \ref{main.thm}. Throughout this section, we assume that $(M,g)$ is a Riemannian manifold satisfying the assumptions of Theorem \ref{main.thm}. 

\begin{lemma}
\label{curvature.inequality}
We have $2 \, R_1(X,JX,Y,JY) \geq -\kappa$ for all points $p \in M$, all $I,J,K \in \mathscr{J}_p$ satisfying $I^2=J^2=K^2=IJK=-\text{\rm id}$, and all unit vectors $X,Y \in T_p M$ satisfying $g(X,Y)=g(IX,Y)=g(JX,Y)=g(KX,Y)=0$. 
\end{lemma} 

\textbf{Proof.} 
Fix a point $p \in M$. Suppose that $I,J,K \in \mathscr{J}_p$ satisfy $I^2=J^2=K^2=IJK=-\text{\rm id}$, and let $X,Y \in T_p M$ be unit vectors satisfying $g(X,Y)=g(IX,Y)=g(JX,Y)=g(KX,Y)=0$. By assumption, $R(X,Y,X,Y) + R(X,JY,X,JY) \geq 0$. Since $R_0(X,Y,X,Y) = R_0(X,JY,X,JY) = \frac{1}{4}$, it follows that 
\[R_1(X,Y,X,Y) + R_1(X,JY,X,JY) \geq -\frac{\kappa}{2}.\] 
On the other hand, using the first Bianchi identity, we obtain 
\begin{align*} 
R_1(X,JX,Y,JY) 
&= -R_1(Y,X,JX,JY) - R_1(JX,Y,X,JY) \\ 
&= R_1(X,Y,X,Y) + R_1(X,JY,X,JY). 
\end{align*} 
Putting these facts together, the assertion follows. This completes the proof of Lemma \ref{curvature.inequality}. \\

\begin{lemma}
\label{four.trace.R_1}
Suppose that $p \in M$, and suppose that $I,J,K \in \mathscr{J}_p$ satisfy $I^2=J^2=K^2=IJK=-\text{\rm id}$. Then 
\begin{align*} 
&R_1(Y,V,W,V) + R_1(Y,IV,W,IV) \\ 
&+ R_1(Y,JV,W,JV) + R_1(Y,KV,W,KV) = 0 
\end{align*}
for all $V,Y,W \in T_p M$.
\end{lemma}

\textbf{Proof.}
Since $R_1$ satisfies the first Bianchi identity, we obtain 
\begin{align*} 
R_1(Y,JW,V,JV) 
&= -R_1(V,Y,JW,JV) - R_1(JW,V,Y,JV) \\ 
&= R_1(Y,V,W,V) + R_1(Y,JV,W,JV). 
\end{align*}
Replacing $V$ by $IV$ gives 
\begin{align*} 
-R_1(Y,JW,IV,KV) 
&= R_1(IV,Y,JW,KV) + R_1(JW,IV,Y,KV) \\ 
&= R_1(Y,IV,W,IV) + R_1(Y,KV,W,KV). 
\end{align*}
Since $R_1(Y,JW,V,JV) = R_1(Y,JW,IV,KV)$, the assertion follows. \\

\begin{lemma}
\label{four.trace.covariant.derivatives.R_1}
Suppose that $p \in M$, and suppose that $I,J,K \in \mathscr{J}_p$ satisfy $I^2=J^2=K^2=IJK=-\text{\rm id}$. Then 
\begin{align*} 
&(D_Z R_1)(Y,V,W,V) + (D_Z R_1)(Y,IV,W,IV) \\ 
&+ (D_Z R_1)(Y,JV,W,JV) + (D_Z R_1)(Y,KV,W,KV) = 0 
\end{align*}
and 
\begin{align*} 
&(D_{X,Z}^2 R_1)(Y,V,W,V) + (D_{X,Z}^2 R_1)(Y,IV,W,IV) \\ 
&+ (D_{X,Z}^2 R_1)(Y,JV,W,JV) + (D_{X,Z}^2 R_1)(Y,KV,W,KV) = 0 
\end{align*}
for all $V,X,Y,Z,W \in T_p M$.
\end{lemma} 

\textbf{Proof.} 
This follows from Lemma \ref{four.trace.R_1}. \\

We next define 
\begin{align*} 
&T_V^{(0)}(X,Y,Z,W) \\ 
&= \sum_{l=1}^{4m} R_0(X,Y,V,e_l) \, R_1(Z,W,V,e_l) \\ 
&+ \sum_{l=1}^{4m} R_0(X,V,Z,e_l) \, R_1(Y,V,W,e_l) - \sum_{l=1}^{4m} R_0(X,V,W,e_l) \, R_1(Y,V,Z,e_l) \\ 
&- \sum_{l=1}^{4m} R_0(Y,V,Z,e_l) \, R_1(X,V,W,e_l) + \sum_{l=1}^{4m} R_0(Y,V,W,e_l) \, R_1(X,V,Z,e_l) 
\end{align*} 
and 
\begin{align*} 
&T_V^{(1)}(X,Y,Z,W) \\ 
&= \sum_{l=1}^{4m} R_1(X,Y,V,e_l) \, R_1(Z,W,V,e_l) \\ 
&+ \sum_{l=1}^{4m} R_1(X,V,Z,e_l) \, R_1(Y,V,W,e_l) - \sum_{l=1}^{4m} R_1(X,V,W,e_l) \, R_1(Y,V,Z,e_l) \\ 
&- \sum_{l=1}^{4m} R_1(Y,V,Z,e_l) \, R_1(X,V,W,e_l) + \sum_{l=1}^{4m} R_1(Y,V,W,e_l) \, R_1(X,V,Z,e_l) 
\end{align*} 
for all points $p \in M$ and all $V,X,Y,Z,W \in T_p M$. For abbreviation, we put 
\[T_V(X,Y,Z,W) = T_V^{(1)}(X,Y,Z,W) + \kappa \, T_V^{(0)}(X,Y,Z,W)\] 
for all points $p \in M$ and all $V,X,Y,Z,W \in T_p M$. 

We further define 
\begin{align*} 
&S_V^{(0)}(X,Y,Z,W) \\ 
&= \sum_{l=1}^{4m} R_0(X,V,e_l,V) \, R_1(e_l,Y,Z,W) + \sum_{l=1}^{4m} R_0(Y,V,e_l,V) \, R_1(X,e_l,Z,W) 
\end{align*}
and 
\begin{align*} 
&S_V^{(1)}(X,Y,Z,W) \\ 
&= \sum_{l=1}^{4m} R_1(X,V,e_l,V) \, R_1(e_l,Y,Z,W) + \sum_{l=1}^{4m} R_1(Y,V,e_l,V) \, R_1(X,e_l,Z,W) 
\end{align*}
for all points $p \in M$ and all $V,X,Y,Z,W \in T_p M$. For abbreviation, we put 
\[S_V(X,Y,Z,W) = S_V^{(1)}(X,Y,Z,W) + \kappa \, S_V^{(0)}(X,Y,Z,W)\]
for all points $p \in M$ and all $V,X,Y,Z,W \in T_p M$. 

\begin{lemma} 
\label{commutator.identity}
We have 
\begin{align*} 
&(D_{X,Z}^2 R_1)(Y,V,W,V) - (D_{X,W}^2 R_1)(Y,V,Z,V) \\ 
&- (D_{Y,Z}^2 R_1)(X,V,W,V) + (D_{Y,W}^2 R_1)(X,V,Z,V) \\ 
&= (D_{V,V}^2 R_1)(X,Y,Z,W) + T_V(X,Y,Z,W) - S_V(X,Y,Z,W) 
\end{align*} 
for all $V,X,Y,Z,W \in T_p M$.
\end{lemma}

\textbf{Proof.} 
We compute 
\begin{align*} 
&(D_{X,V}^2 R_1)(V,Y,Z,W) - (D_{V,X}^2 R_1)(Y,V,W,Z) \\ 
&= \sum_{l=1}^{4m} R(X,V,V,e_l) \, R_1(e_l,Y,Z,W) + \sum_{l=1}^{4m} R(X,V,Y,e_l) \, R_1(V,e_l,Z,W) \\ 
&+ \sum_{l=1}^{4m} R(X,V,Z,e_l) \, R_1(V,Y,e_l,W) + \sum_{l=1}^{4m} R(X,V,W,e_l) \, R_1(V,Y,Z,e_l) 
\end{align*} 
and 
\begin{align*} 
&(D_{Y,V}^2 R_1)(X,V,Z,W) - (D_{V,Y}^2 R_1)(X,V,Z,W) \\ 
&= \sum_{l=1}^{4m} R(Y,V,X,e_l) \, R_1(e_l,V,Z,W) + \sum_{l=1}^{4m} R(Y,V,V,e_l) \, R_1(X,e_l,Z,W) \\ 
&+ \sum_{l=1}^{4m} R(Y,V,Z,e_l) \, R_1(X,V,e_l,W) + \sum_{l=1}^{4m} R(Y,V,W,e_l) \, R_1(X,V,Z,e_l). 
\end{align*} 
Moreover, $R(X,V,Y,e_l) - R(Y,V,X,e_l) = R(X,Y,V,e_l)$ by the first Bianchi identity. Putting these facts together, we obtain 
\begin{align*} 
&(D_{X,V}^2 R_1)(V,Y,Z,W) - (D_{V,X}^2 R_1)(Y,V,W,Z) \\ 
&+ (D_{Y,V}^2 R_1)(X,V,Z,W) - (D_{V,Y}^2 R_1)(X,V,Z,W) \\ 
&= \sum_{l=1}^{4m} R(X,Y,V,e_l) \, R_1(V,e_l,Z,W) \\ 
&+ \sum_{l=1}^{4m} R(X,V,Z,e_l) \, R_1(Y,V,W,e_l) - \sum_{l=1}^{4m} R(X,V,W,e_l) \, R_1(Y,V,Z,e_l) \\ 
&- \sum_{l=1}^{4m} R(Y,V,Z,e_l) \, R_1(X,V,W,e_l) + \sum_{l=1}^{4m} R(Y,V,W,e_l) \, R_1(X,V,Z,e_l) \\ 
&- \sum_{l=1}^{4m} R(X,V,e_l,V) \, R_1(e_l,Y,Z,W) - \sum_{l=1}^{4m} R(Y,V,e_l,V) \, R_1(X,e_l,Z,W) \\ 
&= T_V(X,Y,Z,W) - S_V(X,Y,Z,W). 
\end{align*}
On the other hand, since $R_0$ is parallel, the covariant derivatives of $R_1$ satisfy the second Bianchi identity. This implies 
\[(D_{X,Z}^2 R_1)(Y,V,W,V) - (D_{X,W}^2 R_1)(Y,V,Z,V) = (D_{X,V}^2 R_1)(Y,V,W,Z),\] 
\[-(D_{Y,Z}^2 R_1)(X,V,W,V) + (D_{Y,W}^2 R_1)(X,V,Z,V) = (D_{Y,V}^2 R_1)(X,V,Z,W),\] 
and 
\[(D_{V,V}^2 R_1)(X,Y,Z,W) = (D_{V,X}^2 R_1)(V,Y,Z,W) + (D_{V,Y}^2 R_1)(X,V,Z,W).\] 
From this, the assertion follows. This completes the proof of Lemma \ref{commutator.identity}. \\

\begin{lemma} 
\label{four.trace.commutator.identity}
Suppose that $p \in M$, and suppose that $I,J,K \in \mathscr{J}_p$ satisfy $I^2=J^2=K^2=IJK=-\text{\rm id}$. Then 
\begin{align*} 
0 &= (D_{V,V}^2 R_1)(X,Y,Z,W) + (D_{IV,IV}^2 R_1)(X,Y,Z,W) \\ 
&+ (D_{JV,JV}^2 R_1)(X,Y,Z,W) + (D_{KV,KV}^2 R_1)(X,Y,Z,W) \\ 
&+ T_V(X,Y,Z,W) + T_{IV}(X,Y,Z,W) \\ 
&+ T_{JV}(X,Y,Z,W) + T_{KV}(X,Y,Z,W) \\ 
&- S_V(X,Y,Z,W) - S_{IV}(X,Y,Z,W) \\ 
&- S_{JV}(X,Y,Z,W) - S_{KV}(X,Y,Z,W) 
\end{align*} 
for all $V,X,Y,Z,W \in T_p M$ and all $J \in \mathscr{J}_p$.
\end{lemma} 

\textbf{Proof.} 
This follows by combining Lemma \ref{four.trace.covariant.derivatives.R_1} and Lemma \ref{commutator.identity}. \\

We define 
\[\mu = \sup \{R_1(X,JX,X,JX): \text{\rm $p \in M$, $J \in \mathscr{J}_p$, $X \in T_p M$, $|X|=1$}\}.\] 
The following result follows from work of Berger. \\

\begin{proposition}[cf. M.~Berger \cite{Berger2}]
\label{mu}
We have $\mu=0$ or $\mu=\kappa$.
\end{proposition}

\textbf{Proof.} 
Let us fix a point $p \in M$, an almost complex structure $J \in \mathscr{J}_p$, and a unit vector $X \in T_p M$ such that $R_1(X,JX,X,JX) = \mu$. Moreover, we can find almost complex structures $I,K \in \mathscr{J}_p$ so that $I^2=J^2=K^2=IJK=-\text{\rm id}$. Let $w_1=X$ and $w_2=IX$. Suppose that $w_3,\hdots,w_{2m}$ are vectors in $T_p M$ with the property that $\{w_1,Jw_1,w_2,Jw_2,\hdots,w_{2m},Jw_{2m}\}$ is an orthonormal basis of $T_p M$. In view of the maximality of $R_1(X,JX,X,JX)$, we obtain 
\[\sum_{\alpha=3}^{2m} (D_{w_\alpha,w_\alpha}^2 R_1)(X,JX,X,JX) + \sum_{\alpha=3}^{2m} (D_{Jw_\alpha,Jw_\alpha}^2 R_1)(X,JX,X,JX) \leq 0.\] 
On the other hand, Lemma \ref{four.trace.commutator.identity} implies 
\begin{align*} 
0 
&= \sum_{\alpha=3}^{2m} (D_{w_\alpha,w_\alpha}^2 R_1)(X,JX,X,JX) + \sum_{\alpha=3}^{2m} (D_{Jw_\alpha,Jw_\alpha}^2 R_1)(X,JX,X,JX) \\ 
&+ \sum_{\alpha=3}^{2m} T_{w_\alpha}(X,JX,X,JX) + \sum_{\alpha=3}^{2m} T_{Jw_\alpha}(X,JX,X,JX) \\ 
&- \sum_{\alpha=3}^{2m} S_{w_\alpha}(X,JX,X,JX) - \sum_{\alpha=3}^{2m} S_{Jw_\alpha}(X,JX,X,JX). 
\end{align*}
Putting these facts together, we obtain 
\begin{align} 
\label{inequality.S.T}
&\sum_{\alpha=3}^{2m} S_{w_\alpha}(X,JX,X,JX) + \sum_{\alpha=3}^{2m} S_{Jw_\alpha}(X,JX,X,JX) \notag \\ 
&\leq \sum_{\alpha=3}^{2m} T_{w_\alpha}(X,JX,X,JX) + \sum_{\alpha=3}^{2m} T_{Jw_\alpha}(X,JX,X,JX). 
\end{align} 
Both sides of the inequality (\ref{inequality.S.T}) are independent of the choice of $w_3,\hdots,w_{2m}$. In the following, we choose $w_3,\hdots,w_{2m}$ so that 
\[R_1(X,JX,w_\alpha,w_\beta) = R_1(X,JX,w_\alpha,Jw_\beta) = 0\] 
for all $\alpha,\beta \in \{3,\hdots,2m\}$ with $\alpha \neq \beta$. Applying Lemma 9.16 in \cite{Brendle} to the algebraic curvature tensor $R_1$, we obtain 
\[R_1(X,JX,X,w_\beta) = R_1(X,JX,X,Jw_\beta) = 0\] 
for all $\beta \in \{2,\hdots,2m\}$. Moreover, 
\[R_1(X,JX,X,Iw_\beta) = R_1(X,JX,X,JIw_\beta) = 0\] 
for all $\beta \in \{3,\hdots,2m\}$. Putting these facts together, we conclude that 
\[R_1(X,JX,w_\alpha,w_\beta) = R_1(X,JX,w_\alpha,Jw_\beta) = 0\] 
for all $\alpha,\beta \in \{1,\hdots,2m\}$ with $\alpha \neq \beta$. 

We compute 
\begin{align*} 
&T_{w_\alpha}^{(1)}(X,JX,X,JX) + T_{Jw_\alpha}^{(1)}(X,JX,X,JX) \\ 
&= 2 \sum_{\beta=1}^{2m} R_1(X,JX,w_\alpha,w_\beta)^2 + 2 \sum_{\beta=1}^{2m} R_1(X,JX,w_\alpha,Jw_\beta)^2 \\ 
&- 8 \sum_{\beta=1}^{2m} R_1(X,w_\alpha,JX,w_\beta) \, R_1(JX,w_\alpha,X,w_\beta) \\ 
&- 8 \sum_{\beta=1}^{2m} R_1(X,w_\alpha,JX,Jw_\beta) \, R_1(JX,w_\alpha,X,Jw_\beta) 
\end{align*} 
for each $\alpha \in \{3,\hdots,2m\}$. Using the inequality 
\begin{align*} 
&- 8 \sum_{\beta=1}^{2m} R_1(X,w_\alpha,JX,w_\beta) \, R_1(JX,w_\alpha,X,w_\beta) \\ 
&\leq 2 \sum_{\beta=1}^{2m} (R_1(X,w_\alpha,JX,w_\beta) - R_1(JX,w_\alpha,X,w_\beta))^2 \\ 
&= 2 \sum_{\beta=1}^{2m} R_1(X,JX,w_\alpha,w_\beta)^2 
\end{align*} 
and the inequality 
\begin{align*} 
&- 8 \sum_{\beta=1}^{2m} R_1(X,w_\alpha,JX,Jw_\beta) \, R_1(JX,w_\alpha,X,Jw_\beta) \\ 
&\leq 2 \sum_{\beta=1}^{2m} (R_1(X,w_\alpha,JX,Jw_\beta) - R_1(JX,w_\alpha,X,Jw_\beta))^2 \\ 
&= 2 \sum_{\beta=1}^{2m} R_1(X,JX,w_\alpha,Jw_\beta)^2, 
\end{align*} 
we deduce that 
\begin{align} 
\label{T1.individual.terms}
&T_{w_\alpha}^{(1)}(X,JX,X,JX) + T_{Jw_\alpha}^{(1)}(X,JX,X,JX) \notag \\ 
&\leq 4 \sum_{\beta=1}^{2m} R_1(X,JX,w_\alpha,w_\beta)^2 + 4 \sum_{\beta=1}^{2m} R_1(X,JX,w_\alpha,Jw_\beta)^2 
\end{align} 
for each $\alpha \in \{3,\hdots,2m\}$. We sum both sides of the inequality (\ref{T1.individual.terms}) over $\alpha \in \{3,\hdots,2m\}$. We recall that $R_1(X,JX,w_\alpha,w_\beta) = R_1(X,JX,w_\alpha,Jw_\beta) = 0$ for all $\alpha,\beta \in \{1,\hdots,2m\}$ with $\alpha \neq \beta$. This gives 
\begin{align} 
\label{T1}
&\sum_{\alpha=3}^{2m} T_{w_\alpha}^{(1)}(X,JX,X,JX) + \sum_{\alpha=3}^{2m} T_{Jw_\alpha}^{(1)}(X,JX,X,JX) \notag \\ 
&\leq 4 \sum_{\alpha=3}^{2m} R_1(X,JX,w_\alpha,Jw_\alpha)^2. 
\end{align}
Moreover, using Lemma \ref{four.trace.R_1}, we compute 
\begin{equation} 
\label{S1} 
\sum_{\alpha=3}^{2m} S_{w_\alpha}^{(1)}(X,JX,X,JX) + \sum_{\alpha=3}^{2m} S_{Jw_\alpha}^{(1)}(X,JX,X,JX) = 0. 
\end{equation}
In the next step, we recall a result from \cite{Brendle}. It follows from Proposition 9.23 in \cite{Brendle} that $B(R_1,R_0) = 0$, where $B(R_1,R_0)$ denotes the algebraic expression defined on p.~141 of \cite{Brendle}. This implies 
\begin{equation} 
\label{T0.first.term}
\sum_{\alpha=1}^{2m} T_{w_\alpha}^{(0)}(X,JX,X,JX) + \sum_{\alpha=1}^{2m} T_{Jw_\alpha}^{(0)}(X,JX,X,JX) = 0. 
\end{equation}
Using the formula for $R_0$, we compute 
\begin{align*} 
&T_X^{(0)}(X,JX,X,JX) = R_1(X,JX,X,JX), \\ 
&T_{JX}^{(0)}(X,JX,X,JX) = R_1(JX,X,JX,X), \\ 
&T_{IX}^{(0)}(X,JX,X,JX) = R_1(X,IX,X,IX) + R_1(JX,IX,JX,IX), \\ 
&T_{KX}^{(0)}(X,JX,X,JX) = R_1(X,KX,X,KX) + R_1(JX,KX,JX,KX). 
\end{align*} 
We add these four identities. Using Lemma \ref{four.trace.R_1}, we deduce that 
\begin{equation} 
\label{T0.second.term}
\sum_{\alpha=1}^2 T_{w_\alpha}^{(0)}(X,JX,X,JX) + \sum_{\alpha=1}^2 T_{Jw_\alpha}^{(0)}(X,JX,X,JX) = 0. 
\end{equation} 
Subtracting (\ref{T0.second.term}) from (\ref{T0.first.term}) gives 
\begin{equation}
\label{T0}
\sum_{\alpha=3}^{2m} T_{w_\alpha}^{(0)}(X,JX,X,JX) + \sum_{\alpha=3}^{2m} T_{Jw_\alpha}^{(0)}(X,JX,X,JX) = 0.
\end{equation} 
Since the Ricci tensor of $R_0$ is given by $(m+2) \, g$, we obtain 
\begin{align} 
\label{S0.first.term}
&\sum_{\alpha=1}^{2m} S_{w_\alpha}^{(0)}(X,JX,X,JX) + \sum_{\alpha=1}^{2m} S_{Jw_\alpha}^{(0)}(X,JX,X,JX) \notag \\ 
&= (2m+4) \, R_1(X,JX,X,JX) = (2m+4) \mu. 
\end{align}
Moreover, 
\begin{align} 
\label{S0.second.term}
&\sum_{\alpha=1}^2 S_{w_\alpha}^{(0)}(X,JX,X,JX) + \sum_{\alpha=1}^2 S_{Jw_\alpha}^{(0)}(X,JX,X,JX) \notag \\ 
&= 6 \, R_1(X,JX,X,JX) = 6\mu. 
\end{align} 
Subtracting (\ref{S0.second.term}) from (\ref{S0.first.term}) gives 
\begin{equation} 
\label{S0}
\sum_{\alpha=3}^{2m} S_{w_\alpha}^{(0)}(X,JX,X,JX) + \sum_{\alpha=3}^{2m} S_{Jw_\alpha}^{(0)}(X,JX,X,JX) = (2m-2) \mu.
\end{equation} 
Combining (\ref{inequality.S.T}), (\ref{T1}), (\ref{S1}), (\ref{T0}), and (\ref{S0}), we conclude that 
\begin{equation} 
\label{key.inequality}
(2m-2)\kappa\mu \leq 4 \sum_{\alpha=3}^{2m} R_1(X,JX,w_\alpha,Jw_\alpha)^2. 
\end{equation}
On the other hand, it follows from Lemma \ref{curvature.inequality} that 
\[2 \, R_1(X,JX,w_\alpha,Jw_\alpha) \geq -\kappa\] 
and 
\[2 \, R_1(X,JX,Iw_\alpha,JIw_\alpha) \geq -\kappa\] 
for each $\alpha \in \{3,\hdots,2m\}$. Therefore, 
\begin{equation} 
\label{R_1.term.bounded.by.kappa}
2 \, |R_1(X,JX,w_\alpha,Jw_\alpha)| \leq \kappa 
\end{equation}
for each $\alpha \in \{3,\hdots,2m\}$. Applying Lemma 9.16 in \cite{Brendle} to the algebraic curvature tensor $R_1$, we obtain 
\[2 \, R_1(X,JX,w_\alpha,Jw_\alpha) \leq R_1(X,JX,X,JX) = \mu\] 
and 
\[2 \, R_1(X,JX,Iw_\alpha,JIw_\alpha) \leq R_1(X,JX,X,JX) = \mu\] 
for each $\alpha \in \{3,\hdots,2m\}$. Consequently, 
\begin{equation} 
\label{R_1.term.bounded.by.mu}
2 \, |R_1(X,JX,w_\alpha,Jw_\alpha)| \leq \mu 
\end{equation}
for each $\alpha \in \{3,\hdots,2m\}$. Combining (\ref{key.inequality}) and (\ref{R_1.term.bounded.by.kappa}), we obtain $(2m-2)\kappa\mu \leq (2m-2) \kappa^2$. Combining (\ref{key.inequality}) and (\ref{R_1.term.bounded.by.mu}), we obtain $(2m-2)\kappa\mu \leq (2m-2) \mu^2$. Since $m \geq 2$, $\kappa>0$, and $\mu \geq 0$, it follows that either $\mu = 0$ or $\mu=\kappa$. This completes the proof of Proposition \ref{mu}. \\

\begin{remark}
Suppose that $R(X,Y,X,Y) + R(X,JY,X,JY) > 0$ for all points $p \in M$, all $I,J,K \in \mathscr{J}_p$ satisfying $I^2=J^2=K^2=IJK=-\text{\rm id}$, and all unit vectors $X,Y \in T_p M$ satisfying $g(X,Y)=g(IX,Y)=g(JX,Y)=g(KX,Y)=0$. In this case, the strict inequality holds in (\ref{R_1.term.bounded.by.kappa}). The inequality (\ref{key.inequality}) then implies $(2m-2)\kappa\mu < (2m-2)\kappa^2$. On the other hand, combining (\ref{key.inequality}) and (\ref{R_1.term.bounded.by.mu}) gives $(2m-2)\kappa\mu \leq (2m-2) \mu^2$. Since $m \geq 2$, $\kappa>0$, and $\mu \geq 0$, it follows that $\mu=0$ and $R_1$ vanishes identically. This recovers Berger's classical theorem \cite{Berger2}.
\end{remark}

We now continue with the proof of Theorem \ref{main.thm}. We consider the bundle $\text{\rm End}(TM) \oplus TM$. Let $\pi: \text{\rm End}(TM) \oplus TM \to M$ denote the associated bundle projection. We shall define a subset $\mathscr{E}$ of the total space of this bundle. The set $\mathscr{E}$ consists of all triplets $(p,J,X)$ with the property that $p \in M$, $J \in \mathscr{J}_p$, $X \in T_p M$, and $|X|=1$. Clearly, $\mathscr{E}$ is a compact submanifold of dimension $8m+1$. Moreover, the set $\mathscr{E}$ is invariant under parallel transport. 

We consider the function 
\[\varphi: \mathscr{E} \to \mathbb{R}, \quad (p,J,X) \mapsto \kappa - R_1(X,JX,X,JX).\]

\begin{lemma}
\label{pde.for.varphi}
Let $U$ be an open subset of $M$ which is diffeomorphic to $B^{4m}$. Let $\{e_1,\hdots,e_{4m}\}$ be a local orthonormal frame defined on $U$. For each $l \in \{1,\hdots,4m\}$, we lift the vector field $e_l$ on $U$ to a horizontal vector field on $\mathscr{E} \cap \pi^{-1}(U)$, which we denote by $\mathscr{X}_l$. Moreover, for each $l \in \{1,\hdots,4m\}$, we lift the vector field $D_{e_l} e_l$ on $U$ to a horizontal vector field on $\mathscr{E} \cap \pi^{-1}(U)$, which we denote by $\mathscr{Y}_l$. Then 
\begin{align*} 
&\sum_{l=1}^{4m} \mathscr{X}_l(\mathscr{X}_l(\varphi)) - \sum_{l=1}^{4m} \mathscr{Y}_l(\varphi) \\ 
&\leq -C \inf_{\xi \in \mathscr{V}_{(p,J,X)}, \, |\xi| \leq 1} (D^2 \varphi)(\xi,\xi) + C \sup_{\xi \in \mathscr{V}_{(p,J,X)}, \, |\xi| \leq 1} d\varphi(\xi) + C \, |\varphi| 
\end{align*}
at each point $(p,J,X) \in \mathscr{E} \cap \pi^{-1}(U)$. Here, $\mathscr{V}_{(p,J,X)} \subset T_{(p,J,X)} \mathscr{E}$ denotes the vertical space at the point $(p,J,X) \in \mathscr{E}$. Moreover, $C$ is a large positive constant. 
\end{lemma}

\textbf{Proof.} 
The Ricci tensor of $(M,g)$ is given by $(m+2)\kappa \, g$. Using Proposition 2.11 in \cite{Brendle}, we obtain 
\[\Delta R + Q(R) = (2m+4)\kappa \, R.\] 
Here, 
\[\Delta R = \sum_{l=1}^{4m} D_{e_l,e_l}^2 R = \sum_{l=1}^{4m} D_{e_l} D_{e_l} R - \sum_{l=1}^{4m} D_{D_{e_l} e_l} R\] 
denotes the Laplacian of the Riemann curvature tensor, and $Q(R)$ denotes the algebraic expression defined on p.~53 of \cite{Brendle}. Clearly, $\Delta R = \Delta R_1$ since $R_0$ is parallel. Moreover, Proposition 9.23 in \cite{Brendle} implies that $Q(R) = Q(R_1) + \kappa^2 \, Q(R_0) = Q(R_1) + (2m+4) \kappa^2 \, R_0$. Putting these facts together, we conclude that 
\begin{equation} 
\label{pde.for.R_1}
\Delta R_1 + Q(R_1) = (2m+4)\kappa \, R_1. 
\end{equation}
For each $l \in \{1,\hdots,4m\}$, we have 
\[\mathscr{X}_l(\mathscr{X}_l(\varphi)) = -(D_{e_l} D_{e_l} R_1)(X,JX,X,JX)\] 
and 
\[\mathscr{Y}_l(\varphi) = -(D_{D_{e_l} e_l} R_1)(X,JX,X,JX)\] 
at each point $(p,J,X) \in \mathscr{E} \cap \pi^{-1}(U)$. This implies 
\begin{equation} 
\label{horizontal.Laplacian.of.varphi.1}
\sum_{l=1}^{4m} \mathscr{X}_l(\mathscr{X}_l(\varphi)) - \sum_{l=1}^{4m} \mathscr{Y}_l(\varphi) = -\Delta R_1(X,JX,X,JX) 
\end{equation} 
at each point $(p,J,X) \in \mathscr{E} \cap \pi^{-1}(U)$. Combining (\ref{pde.for.R_1}) and (\ref{horizontal.Laplacian.of.varphi.1}), we obtain 
\begin{align} 
\label{horizontal.Laplacian.of.varphi.2}
&\sum_{l=1}^{4m} \mathscr{X}_l(\mathscr{X}_l(\varphi)) - \sum_{l=1}^{4m} \mathscr{Y}_l(\varphi) \notag \\ 
&= Q(R_1)(X,JX,X,JX) - (2m+4)\kappa \, R_1(X,JX,X,JX) 
\end{align}
at each point $(p,J,X) \in \mathscr{E} \cap \pi^{-1}(U)$. 

In order to estimate the term on the right hand side of (\ref{horizontal.Laplacian.of.varphi.2}), we adapt the proof of Theorem 9.24 in \cite{Brendle}. Let us fix a triplet $(p,J,X) \in \mathscr{E} \cap \pi^{-1}(U)$. As above, we put $w_1=X$ and $w_2=IX$. We can find vectors $w_3,\hdots,w_{2m} \in T_p M$ with the property that $\{w_1,Jw_1,w_2,Jw_2,\hdots,w_{2m},Jw_{2m}\}$ is an orthonormal basis of $T_p M$ and 
\begin{equation} 
\label{diagonal}
R_1(X,JX,w_\alpha,w_\beta) = R_1(X,JX,w_\alpha,Jw_\beta) = 0 
\end{equation}
for all $\alpha,\beta \in \{3,\hdots,2m\}$ with $\alpha \neq \beta$. 

By adapting the proof of Lemma 9.16 in \cite{Brendle}, we can show that 
\begin{equation} 
\label{first.variation.a}
|R_1(X,JX,X,w_\beta)| \leq C \sup_{\xi \in \mathscr{V}_{(p,J,X)}, \, |\xi| \leq 1} d\varphi(\xi) 
\end{equation}
and 
\begin{equation} 
\label{first.variation.b}
|R_1(X,JX,X,Jw_\beta)| \leq C \sup_{\xi \in \mathscr{V}_{(p,J,X)}, \, |\xi| \leq 1} d\varphi(\xi) 
\end{equation}
for all $\beta \in \{2,\hdots,2m\}$. An analogous argument gives 
\[|R_1(X,JX,X,Iw_\beta)| \leq C \sup_{\xi \in \mathscr{V}_{(p,J,X)}, \, |\xi| \leq 1} d\varphi(\xi)\] 
and 
\[|R_1(X,JX,X,JIw_\beta)| \leq C \sup_{\xi \in \mathscr{V}_{(p,J,X)}, \, |\xi| \leq 1} d\varphi(\xi)\] 
for all $\beta \in \{3,\hdots,2m\}$. These inequalities can be rewritten as 
\begin{equation} 
\label{first.variation.c}
|R_1(X,JX,IX,w_\beta)| \leq C \sup_{\xi \in \mathscr{V}_{(p,J,X)}, \, |\xi| \leq 1} d\varphi(\xi) 
\end{equation}
and 
\begin{equation} 
\label{first.variation.d}
|R_1(X,JX,IX,Jw_\beta)| \leq C \sup_{\xi \in \mathscr{V}_{(p,J,X)}, \, |\xi| \leq 1} d\varphi(\xi) 
\end{equation}
for all $\beta \in \{3,\hdots,2m\}$. Combining (\ref{diagonal}), (\ref{first.variation.a}), (\ref{first.variation.b}), (\ref{first.variation.c}), and (\ref{first.variation.d}), we conclude that 
\begin{equation} 
\label{estimate.for.R1.a}
|R_1(X,JX,w_\alpha,w_\beta)| \leq C \sup_{\xi \in \mathscr{V}_{(p,J,X)}, \, |\xi| \leq 1} d\varphi(\xi) 
\end{equation}
and 
\begin{equation} 
\label{estimate.for.R1.b}
|R_1(X,JX,w_\alpha,Jw_\beta)| \leq C \sup_{\xi \in \mathscr{V}_{(p,J,X)}, \, |\xi| \leq 1} d\varphi(\xi) 
\end{equation}
for all $\alpha,\beta \in \{1,\hdots,2m\}$ with $\alpha \neq \beta$. Moreover, by adapting the proof of Lemma 9.16 in \cite{Brendle}, we can show that  
\begin{align} 
\label{second.variation.a}
2 \, R_1(X,JX,w_\alpha,Jw_\alpha) 
&\leq R_1(X,JX,X,JX) \notag \\ 
&- C \inf_{\xi \in \mathscr{V}_{(p,J,X)}, \, |\xi| \leq 1} (D^2 \varphi)(\xi,\xi) \\ 
&+ C \sup_{\xi \in \mathscr{V}_{(p,J,X)}, \, |\xi| \leq 1} d\varphi(\xi) \notag
\end{align}
for each $\alpha \in \{3,\hdots,2m\}$. An analogous argument gives 
\begin{align*} 
2 \, R_1(X,JX,Iw_\alpha,JIw_\alpha) 
&\leq R_1(X,JX,X,JX) \\ 
&- C \inf_{\xi \in \mathscr{V}_{(p,J,X)}, \, |\xi| \leq 1} (D^2 \varphi)(\xi,\xi) \\ 
&+ C \sup_{\xi \in \mathscr{V}_{(p,J,X)}, \, |\xi| \leq 1} d\varphi(\xi) 
\end{align*}
for each $\alpha \in \{3,\hdots,2m\}$. This inequality can be rewritten as 
\begin{align} 
\label{second.variation.b}
-2 \, R_1(X,JX,w_\alpha,Jw_\alpha) 
&\leq R_1(X,JX,X,JX) \notag \\ 
&- C \inf_{\xi \in \mathscr{V}_{(p,J,X)}, \, |\xi| \leq 1} (D^2 \varphi)(\xi,\xi) \\
&+ C \sup_{\xi \in \mathscr{V}_{(p,J,X)}, \, |\xi| \leq 1} d\varphi(\xi) \notag
\end{align}
for each $\alpha \in \{3,\hdots,2m\}$. Combining (\ref{second.variation.a}) and (\ref{second.variation.b}), we conclude that 
\begin{align} 
\label{estimate.for.R1.c}
2 \, |R_1(X,JX,w_\alpha,Jw_\alpha)| 
&\leq R_1(X,JX,X,JX) \notag \\ 
&- C \inf_{\xi \in \mathscr{V}_{(p,J,X)}, \, |\xi| \leq 1} (D^2 \varphi)(\xi,\xi) \\ 
&+ C \sup_{\xi \in \mathscr{V}_{(p,J,X)}, \, |\xi| \leq 1} d\varphi(\xi) \notag
\end{align}
for each $\alpha \in \{3,\hdots,2m\}$. 

Proposition 9.15 in \cite{Brendle} gives 
\begin{align*} 
Q(R_1)(X,JX,X,JX) 
&\leq -2 \, R_1(X,JX,X,JX)^2 \\ 
&+ 4 \sum_{\alpha,\beta=1}^{2m} R_1(X,JX,w_\alpha,w_\beta)^2 \\ 
&+ 4 \sum_{\alpha,\beta=1}^{2m} R_1(X,JX,w_\alpha,Jw_\beta)^2. 
\end{align*}
Using (\ref{estimate.for.R1.a}) and (\ref{estimate.for.R1.b}), we obtain 
\begin{align*} 
Q(R_1)(X,JX,X,JX) 
&\leq 6 \, R_1(X,JX,X,JX)^2 \\ 
&+ 4 \sum_{\alpha=3}^{2m} R_1(X,JX,w_\alpha,Jw_\alpha)^2 \\ 
&+ C \sup_{\xi \in \mathscr{V}_{(p,J,X)}, \, |\xi| \leq 1} d\varphi(\xi). 
\end{align*}
Using (\ref{estimate.for.R1.c}), we conclude that 
\begin{align} 
\label{estimate.for.Q}
Q(R_1)(X,JX,X,JX) 
&\leq (2m+4) \, R_1(X,JX,X,JX)^2 \notag \\ 
&- C \inf_{\xi \in \mathscr{V}_{(p,J,X)}, \, |\xi| \leq 1} (D^2 \varphi)(\xi,\xi) \\ 
&+ C \sup_{\xi \in \mathscr{V}_{(p,J,X)}, \, |\xi| \leq 1} d\varphi(\xi). \notag
\end{align}
Finally, 
\begin{equation} 
\label{term.on.rhs}
(2m+4) \, R_1(X,JX,X,JX)^2 \leq (2m+4)\kappa \, R_1(X,JX,X,JX) + C \, |\varphi| 
\end{equation}
by definition of $\varphi$. Combining (\ref{horizontal.Laplacian.of.varphi.2}), (\ref{estimate.for.Q}), and (\ref{term.on.rhs}), the assertion follows. \\

\begin{proposition}
\label{invariance.under.parallel.transport}
The function $\varphi: \mathscr{E} \to \mathbb{R}$ is nonnegative. Moreover, the set $\mathscr{F} = \{\varphi=0\} \subset \mathscr{E}$ is invariant under parallel transport.
\end{proposition}

\textbf{Proof.} 
It follows from Proposition \ref{mu} that $\mu \leq \kappa$. This implies $\inf_{\mathscr{E}} \varphi = \kappa-\mu \geq 0$. This proves the first statement. To prove the second statement, suppose that $U$ is an open subset of $M$ which is diffeomorphic to $B^{4m}$. Using Lemma \ref{pde.for.varphi} and the strict maximum principle (cf. \cite{Brendle}, Corollary 9.7), we conclude that every smooth horizontal curve in $\mathscr{E} \cap \pi^{-1}(U)$ that starts at a point in $\mathscr{F} \cap \pi^{-1}(U)$ must be contained in $\mathscr{F} \cap \pi^{-1}(U)$. From this, we deduce that every smooth horizontal curve in $\mathscr{E}$ that starts at a point in $\mathscr{F}$ must be contained in $\mathscr{F}$. This completes the proof of Proposition \ref{invariance.under.parallel.transport}. \\

\begin{proposition}
\label{mu.vanishes}
We have $\mu=0$.
\end{proposition}

\textbf{Proof.} 
We argue by contradiction. Suppose that $\mu \neq 0$. Using Proposition \ref{mu}, we obtain $\mu=\kappa$. Consequently, we can find a point $p \in M$, an almost complex structure $J \in \mathscr{J}_p$, and a unit vector $X \in T_p M$ such that $R_1(X,JX,X,JX) = \kappa$. Moreover, we can find almost complex structures $I,K \in \mathscr{J}_p$ so that $I^2=J^2=K^2=IJK=-\text{\rm id}$. 

Let $Z$ be an arbitrary unit vector in $T_p M$. Since $\text{\rm Sp}(m)$ acts transitively on the unit sphere $S^{4m-1}$, we can find a linear isometry $L: T_p M \to T_p M$ which commutes with $I,J,K$ and satisfies $LX = Z$. Since $\text{\rm Hol}(M,g) = \text{\rm Sp}(m) \cdot \text{\rm Sp}(1)$, there exists a piecewise smooth path $\gamma: [0,1] \to M$ such that $\gamma(0) = \gamma(1) = p$ and $P_\gamma = L$, where $P_\gamma: T_{\gamma(0)} M \to T_{\gamma(1)} M$ denotes the parallel transport along $\gamma$. By 
Proposition \ref{invariance.under.parallel.transport}, the set $\mathscr{F}$ is is invariant under parallel transport. Since $(p,J,X) \in \mathscr{F}$, it follows that $(p,LJL^{-1},LX) \in \mathscr{F}$. Thus, $(p,J,Z) \in \mathscr{F}$. In other words, $R_1(Z,JZ,Z,JZ) = \kappa$. 

To summarize, we have shown that $R_1(Z,JZ,Z,JZ) = \kappa$ for every unit vector $Z \in T_p M$. Since $\kappa>0$, it follows that the scalar curvature of $R_1$ at the point $p$ is strictly positive. This contradicts the fact that the scalar curvature of $R_1$ vanishes identically. This completes the proof of Proposition \ref{mu.vanishes}. \\

Using Proposition \ref{mu.vanishes}, it is easy to see that $R_1=0$ at each point in $M$. This completes the proof of Theorem \ref{main.thm}.

\section{Proof of Corollary \ref{consequence.of.main.thm}}

In this final section, we give the proof of Corollary \ref{consequence.of.main.thm}. Suppose that $(M,g)$ is a Riemannian manifold satisfying the assumptions of Corollary \ref{consequence.of.main.thm}. Since $(M,g)$ is a quaternionic-K\"ahler manifold with positive scalar curvature, $M$ is simply connected (see \cite{Salamon}, Theorem 6.6). As above, we denote by $\text{\rm Hol}(M,g) \subset \text{\rm Sp}(m) \cdot \text{\rm Sp}(1)$ the holonomy group of $(M,g)$. We distinguish two cases: 

\textit{Case 1:} Suppose that $\text{\rm Hol}(M,g) = \text{\rm Sp}(m) \cdot \text{\rm Sp}(1)$. In this case, Theorem \ref{main.thm} implies that $R_1=0$ at each point in $M$. Thus, the Riemann curvature tensor of $(M,g)$ is parallel, and $(M,g)$ is isometric to a symmetric space.

\textit{Case 2:} Suppose that $\text{\rm Hol}(M,g) \subsetneq \text{\rm Sp}(m) \cdot \text{\rm Sp}(1)$. In this case, we obtain $\dim \text{\rm Hol}(M,g) < \dim(\text{\rm Sp}(m) \cdot \text{\rm Sp}(1))$. In particular, $\dim \text{\rm Hol}(M,g) < \dim \text{\rm U}(2m)$ and $\dim \text{\rm Hol}(M,g) < \dim \text{\rm SO}(4m)$. Since $(M,g)$ is a quaternionic-K\"ahler manifold with positive scalar curvature, it follows from a theorem of Berger that $(M,g)$ is irreducible (see \cite{Berger1} or \cite{Besse}, Theorem 14.45(b)). Using Berger's holonomy classification theorem (cf. \cite{Besse}, Corollary 10.92), we conclude that $(M,g)$ is isometric to a symmetric space. This completes the proof of Corollary \ref{consequence.of.main.thm}.

\end{document}